\documentclass[12pt]{article}
\usepackage[dvips]{graphicx}
\usepackage{amssymb}
\usepackage{amsmath,amsthm}
\usepackage{natbib}

\usepackage{bm}
\bmdefine{\bi}{i}
\bmdefine{\bj}{j}

\date{March, 2006}

\title{Evaluation of per-record identification risk and
  swappability of records in a microdata set via decomposable models}
\author{Akimichi Takemura and Yushi Endo\\
  Graduate School of Information Science and Technology\\
  University of Tokyo }
\begin{document}
\maketitle
\begin{abstract}
  We propose a strategy for disclosure risk evaluation and disclosure
  control of a microdata set based on fitting decomposable models of a
  multiway contingency table corresponding to the microdata set.  By
  fitting decomposable models, we can evaluate per-record
  identification (or re-identification) risk of a microdata set.
  Furthermore we can easily determine swappability of risky records
  which does not disturb the set of marginals of the decomposable
  model.  Use of decomposable models has been already considered in
  the existing literature.  The contribution of this paper is to
  propose a systematic strategy to the problem of finding a model with
  a good fit, identifying risky records under the model, and then
  applying the swapping procedure to these records.
\end{abstract}

\section{Introduction}

In this paper we propose a systematic strategy of per-record
identification risk and disclosure control of risky records of a
microdata set by fitting decomposable models to a multiway contingency
tables corresponding to the microdata.  The first stage of our
strategy consists of selecting decomposable models with a good fit to
the data based on Akaike's information criterion (AIC).  Since the
number of decomposable models is large, we propose an algorithm to
find locally optimum decomposable models.  The second stage is to
evaluate cell probabilities of sample unique records and to estimate
the number of population uniques in the microdata set based on the
chosen model.  The third stage consists of disclosure control of risky
records by swapping.  We consider swapping which does not disturb the
set of marginals corresponding to the chosen model.

In evaluating the disclosure risk of a given microdata set, the number
of the population uniques among the sample unique records has been
considered to an important overall measure of the disclosure risk.
Starting from Poisson-Gamma model (\cite{beth:kell:pann:1990}) various
models of random partitions have been proposed for estimating the
number of population uniques.  See a series of works of Hoshino
(\cite{hoshino-takemura}, \cite{hoshino-2001jos},
\cite{hoshino-2003jjss}, \cite{hoshino-2005aism}) and references
therein.  These models treat the sample unique records exchangeably
and hence the estimated conditional probability of population
uniqueness is common for every sample unique record.  However some
sample unique records are clearly more likely to be population uniques
than other records, according to ``rareness'' of the 
records.  If a sample unique has outlying observations or has very
a rare combination of observed characteristics, it is likely to be a
population unique.  A simple descriptive method for evaluating
per-record identification risk is to look at minimum unsafe combination of
variables for a sample unique record (\cite{takemura-unsafe}).

More systematic way of evaluating the per-record identification risk
is to model cell probabilities of the contingency table corresponding
to a microdata set, where all the key variables of the microdata set
are categorized and the joint frequencies of the key variables are
counted.  If the estimated cell probability of a sample unique cell
is very small, then the sample unique is rare and risky.  This approach was
investigated in \cite{skinner-holmes98}, \cite{fienberg-makov98},
\cite{elamir}.  They used the standard log-linear models for cell
probabilities of contingency tables.
% In this paper we consider fitting the Lancaster-type additive model of
% interaction terms, because of its simple computation and
% interpretation.

In actual evaluation of disclosure risk, we often have to consider 10
or more possible key variables.  Then the contingency table is
large and sparse and the 
% In the application in Section 4
% the total number of the cells
% is about 1/4 billion.
estimation of cell probabilities of  standard log-linear models is not
straightforward, except for decomposable models.
In Section \ref{subsec:example} we consider an example
of a 8-way contingency table from 1990 U.S.\ Census of Population and
Housing data.  From the viewpoint of disclosure control this example is 
of moderate size  but the contingency table corresponding to
the microdata has more than 12 million cells.

Because of the computational difficulty Takemura
\cite{takemura02-iass} considered Lancaster-type additive modeling of
cell probabilities.  However in fitting additive models estimated cell
probabilities often become negative, especially for empty cells.  In
this sense additive models are not satisfactory for estimating small
cell probabilities, although they are useful for the purpose of
relative evaluation of identification risks of sample unique cells.

% Throughout this paper we assume the following simple superpopulation
% model: the cell probabilities of the contingency cell are unknown but
% fixed  and each of the $N$ individuals  of the population falls into a
% cell
% by an independent multivariate Bernoulli trial.
% Concerning the sampling we assume simple random sampling of $n$
% individuals without replacement.
% In this setting, the unobserved $N-n$ individuals are distributed
% independently
% of the observed $n$ individuals and the evaluation of conditional
% probability
% is simply derived from the multinomial probability of the unobserved
% individuals.

Among the log-linear models, decomposable models are special in the
sense that the maximum likelihood estimates of the cell probabilities
can be explicitly written as  ratios of products of marginal
frequencies.  Unlike other log-linear models, in a decomposable model
cell probability of each cell can be separately estimated.  This is a
very attractive feature of decomposable model, because we are mainly
interested in sample unique cells or other cells of small
frequency.  Furthermore model selection among decomposable models is
relatively easy, because the maximized log likelihood and the degrees
of freedom can be simply evaluated.  For fitting other log-linear
models, we need some iterative procedure such as iterative
proportional scaling (see e.g.\ \cite{endo-takemura}).  For large
contingency tables iterative proportional scaling is computationally
very intensive, because cell probability estimates of all the cells have
to be stored in some form and updated in each iteration.

Estimation and diagnostics of a particular decomposable model is easy.
However if the number $m$ of key variables is large, there are many
possible decomposable models.  In Table \ref{result1} below, for our
example of $m=8$ key variables, there are more than 30 million
possible decomposable models.  Finding the best fitting model among
more than 30 million possible models is impractical.  We propose to
find several locally optimum models and choose one of these models.

Once a decomposable model with a good fit is obtained, we look at
sample unique cells with very small estimated cell probabilities.  If
the cells are considered to be risky, it is desirable to perform some
disclosure control measure to these cells.  From the viewpoint of
log-linear model, it is natural to consider swapping of these risky
records in such a way that the swapping does not disturb the given set
of marginals corresponding to the cliques of the decomposable model.
This is based on the fact that the set of marginals constitutes the
sufficient statistic of the model and swapping does not influence
statistical inferences based on the model.  Using the results of
\cite{takemura-hara} we show that it is straightforward to determine
whether a particular record is swappable and find another record for
swapping if swapping is possible.

The organization of the paper is as follows. In Section 
\ref{sec:prelimininaries} we summarize preliminary material and
introduce our working example. 
In Section \ref{sec:selection} we discuss fitting and selection of
decomposable models.  In Section 
\ref{sec:risk} based on a chosen decomposable model we evaluate
per-record identification risk.  In  Section
\ref{sec:swappability} we perform  swapping of risky records.
Section \ref{sec:remarks} ends the paper with some concluding remarks.

\section{Preliminaries and a working example}
\label{sec:prelimininaries}

In this section we prepare notations on decomposable models and
describe a working example analyzed in this paper.

\subsection{Notations on decomposable models}
\label{subsec:notation}

We follow the notation of \cite{lauritzen1996}.
Let $\Delta=\{1,\ldots,m\}$ denote the set of the key variables. Each
variable is denoted by $\delta\in \Delta$.    
We assume that all key variables are already discretized and let 
${\cal
  I}_\delta=\{1,\ldots,I_\delta\}$ 
denote the set of categories of $\delta$.
Each cell is indexed by  $m$ indices $\bi=(i_1,\ldots,i_m)$ and 
the set of the cells is the direct product
$
{\cal I} = \prod_{\delta \in \Delta} {\cal I}_\delta
$.
The frequency of cell $\bi$  is denoted by $n(\bi)$.

Let $a \subset \Delta$  be a subset  of variables.
Then an $a$-marginal cell $\bi_a$  of $\bi=(i_1, \dots, i_m) $ is defined  as
${\bi_a} = (i_\delta)_{\delta \in a}$.
The set of $a$-marginal cells is ${\cal I}_a = \prod_{\delta \in a} {\cal I}_\delta$.
The marginal  frequency of 
$a$-marginal cell $\bi_a$ is written as
\[
n(\bi_a)= \sum_{\bj: \bj_a = \bi_a} n(\bj),
\]
where
% $\Delta=\{1,\ldots,m\}$,
% $\bi=(i_1,\ldots,i_m)$, $\bj=(j_1,\ldots,j_m)$ $B$N>l9g$K$O(B
$\bj_a = \bi_a$ means $i_k=j_k, \ \forall k\in a$.
Let $n=\sum_{\bi \in {\cal I}} n(\bi)$ denote the sample size (number
of records) of the
microdata set.  We denote the relative frequency of a cell $\bi$ and a
marginal cell $\bi_a$ by
\[
r(\bi) =\frac{n(\bi)}{n}, \quad  r(\bi_a) =\frac{n(\bi_a)}{n}.
\]
We use the same notation for cell probabilities
$p(\bi)$, $p(\bi_a)$, etc.

Consider a graph $G=(\Delta,E)$ with the set of vertices $\Delta$ and the set of
edges $E$.  Let $\cal C$ denote the set of (maximal) cliques.
For a subset $a\subset\Delta$ let $\mu_a : {\cal I}\rightarrow R$
denote a function of $\bi$ which only depends on the marginal cell
$\bi_a$, i.e.\ $\mu_a(\bi)=\mu_a(\bi_a)$.
Then the graphical model associated with  $G$ specifies  the cell
probability $p(\bi)$  as
\begin{equation}
\label{eq:hierarchical-model}
\log p(\bi) = \sum_{a\in {\cal C}} \mu_a(\bi_a).
\end{equation}

A graph $G$ is {\em chordal} ({\em decomposable}, {\em triangulated}),
if every cycle of length $l\geq 4$ has
a chord.   A graphical model with a chordal $G$ is called a
{\it decomposable} model.
For a decomposable model, the cliques 
can be ordered to satisfy
the running intersection property:
\begin{center}
(RIP) \quad For each $2\leq j\leq m$, there exists $1\leq k\leq j-1$,
such
that $
S_j = C_j\cap (C_1\cup C_2\cup \cdots \cup C_{j-1})\subset C_k.
$
\end{center}
An ordering $(C_1,\dots,C_m)$ satisfying RIP is called a {\it perfect
sequence}.
% For a perfect sequence of a decomposable model let
% \[
% S_j = C_j\cap (C_1\cup C_2\cup \cdots \cup C_{j-1}), \qquad 2\leq
% j\leq m.
% \]
$S_2,\dots,S_m$ are minimal vertex separators of $G$. 
The number of times a minimal vertex
separator $S$ appears in any perfect sequence is the same and called
the {\it multiplicity} of $S$. We denote the multiplicity of $S$ by
$\nu (S)$.  $\cal S$ denotes the set of minimal vertex separators.
% \[
% {\cal S}=\{ S_2, \dots, S_m\},
% \]
% denotes the multiset of minimal vertex separators, where each minimal
% vertex separator $S$ appears $\nu(S)$ times in $\cal S$.
In the following we simply say ``separator'' to mean a minimal vertex
separator.

The maximum likelihood estimate (MLE) of 
a decomposable model is explicitly written
as
% Then MLE of a decomposable
% model is obtained by following formula,
\begin{equation}
\label{eq:mle}
\hat p_{{\rm MLE}}(i)
=\begin{cases}
\displaystyle \frac{\prod _{C\in {\cal C}}r(i_C)}{\prod _{S\in {\cal
      S}}r(i_S)^{\nu(S)}}, &
\text{if}\ \   r(i_C) >0,\ \forall C\in \cal C,\\
0, & \text{otherwise}.
\end{cases}
\end{equation}
The degrees of freedom is also simply written (Proposition 4.35 of
\cite{lauritzen1996}).
\begin{equation}
\label{eq:df}
\sum_{C \in {\cal C}}  \prod_{\delta \in C} I_\delta -
\sum_{S \in {\cal S}} \nu(S) \prod_{\delta \in S} I_\delta .
\end{equation}
Hence AIC  for model selection is
also easily computed.
\begin{equation}
\label{eq:AIC}
\text{AIC}=  -2 \times \text{(log likelihood)} + 2  \times
\text{(degrees of freedom)}.
\end{equation}

\begin{table}[htb]
\caption{Number of graphical models and decomposable models}
\label{tbl:decomposable}
\begin{center}
\begin{tabular}{| c | c| c |}\hline
    $m$ & graphical & decomposable\\ \hline
  2 & 2 & 2  (2) \\ \hline
  3 & 8 &  8  (4)\\ \hline
  4 & 64 & 61 (10) \\ \hline
  5 & 1024 & 820 (27) \\ \hline
  6 & 32768  & 18154 (96)\\ \hline
  7 & 2097152\  & 617675 (469) \\ \hline
  8 & 268435456 & 30888596 (3734) \\ \hline
\end{tabular} 
\end{center} 
\end{table}

In Table \ref{tbl:decomposable} we list the number of graphical models
and the number of decomposable models for $m$-way contingency tables
up to $m=8$.  We see that the number of decomposable models
increases very fast with $m$.  The number in the parentheses for the
decomposable model indicates the number of chordal graphs of $m$
vertices after identification of isomorphic graphs, i.e., we do not
distinguish graphs which can be obtained by relabeling of vertices.
Based on \cite{endo-bachelor-thesis} we provide a list of
non-isomorphic chordal graphs for $m\le 8$ in
\cite{endo-takemura-list}.  Given a list of non-isomorphic chordal
graphs we can pick a decomposable model by choosing an graph from
the list and arbitrary assigning a variable to each vertex of the
graph.

\subsection{A working example}
\label{subsec:example}

In this paper we apply our strategy to a test data set from
1990 U.S.\ Census of Population and Housing Public Use Microdata
Samples.  We subsampled $n=9809$ individuals from the state of
Washington and chose $m=8$ variables for our experiment.
\begin{center}
\begin{tabular}{|l|l|}\hline
1.\ Relationship (14 categories) \qquad &  2.\ Sex (2 categories)\\ \hline
3.\ Age (91 categories) &  4.\  Marital status (5 categories) \\ \hline
5.\ Place of birth (14 categories)\ &  6.\ Spouse present/absent (7
categories) \\ \hline
7.\  Own child (2 categories) & 8.\ Age of own child (5 categories)\\ \hline
% \   9.\ Related child (2),\
% 10.\ Detailed relationship (10).
\end{tabular}
\end{center}
% Type here the text
The population size of the state of Washington is about 
$N={\rm 4,867,000}$. The dataset can be viewed
as a $8$-way  contingency table of the type
%14\times 2 \times 91\times 5 \times 14\times 7 \times 2 \times 5
%\times 2 \times 10
\[
14\times 2 \times 91\times 5 \times
14\times 7 \times 2 \times 5 
\]
with approximately 12.5 million cells (more exactly 12,485,200 cells).
%with 249,704,000 cells.    
We see that the contingency table is very
sparse with only $n=9809$  counts 
in 12.5 million cells.
%among 249,704,000 cells
We took these $m=8$ variables from a PUMS data set without further
global recoding.  For example we used the age itself with 91
categories. This is somewhat unrealistic for evaluation of
disclosure risk.  On the other hand there are other possible key
variables in the original PUMS data set.  
% We also intended to check
% how the proposed model works for large contingency tables.

It should be noted that although the (formal) total number of cells 
12,485,200
% 249,704,000 
is very large, the effective total number should be much smaller because of
structural zeros.  For example there is no age of own child if there
is no own child. In this case the age of own child is coded as N/A in
the original data set.  Also there
is an obvious relation between age and marital status.
% The existence of large number of structural zeros seems to adversely
% affect the fit of the model as discussed below.
In this paper we ignore the effect of structural zeros.
See Section \ref{sec:remarks} for more discussion.

For reference we show first few lines of
%$n=9809{\rm rows} \times m=10{\rm colums}$ data matrix.
$9809 \times 8$ data matrix.
\begin{center}
{\tt
00,0,17,4,10,6,0,0\\
00,0,17,4,52,6,0,0\\
00,0,18,0,23,1,0,0\\
00,0,18,0,24,1,0,0\\
00,0,18,0,51,1,0,0}
\end{center}

The frequencies of the cell sizes
(size indices, frequency of frequencies) % $s_1, s_2, \ldots$, 
of this data set is given as follows.  The table shows that 
there are 2243 cells of frequency 1, 524 cells of frequency 2, etc.
\vspace{5mm}
\begin{center}
\begin{tabular}{|c|ccccccccccc|} \hline
Cell size & 1 & 2   & 3 & 4 & 5 & 6 & 7 & 8 & 9 & 10 & 11 $\le$\\
\hline
Frequency &    2243 & 524 &275&132&104 & 60& 59& 34& 46& 19 & 124 \\ \hline
\end{tabular}
\end{center}
\vspace{5mm}

We are interested in estimating the number of population uniques among
2243 sample uniques and evaluate which sample record is particularly risky.
As a preliminary analysis, we 
fitted   Ewens model, Pitman model and Lancaster-type additive model.
The estimates of the number of population uniques of these models are as follows.
\begin{center}
Ewens model: 5.9, \quad
Pitman model: 214.0, \quad
additive model: 252.1.
\end{center}

\section{Selection of decomposable models}
\label{sec:selection}

The first step of our strategy is to choose a decomposable model which
fits the data.  As shown in Table \ref{tbl:decomposable} the number of possible
decomposable models grow very fast as the number of variables $m$
increases.  For $m\le 8$ we can use the list of non-isomorphic chordal
graphs available at \cite{endo-takemura-list}.  We present the following
Algorithm 1 to obtain locally best decomposable model in terms of AIC.
Application of Algorithm 1 to the data set of our working example is summarized
in Table \ref{result1} below.

In our algorithm we add or subtract an edge to (or from) a chordal
graph to move to another chordal graph and evaluate AIC.  It outputs a
model with locally minimum AIC.  We  can apply our algorithm from
various initial models and compare these locally best models to obtain
approximately a globally best model.

Notations of Algorithm 1 is as follows. 
$G=G(V,E)=G(V,E_G)$ is a graph
with the set of vertices $V$ and the set of vertices $E$.
$M(G)$ denotes the graphical model associated with $G$.
$E(C_m)$ denotes the set of edges of the complete graph with $m$
vertices.

In Step 1 we choose an initial model randomly from the list of
non-isomorphic decomposable models(\cite{endo-bachelor-thesis},
\cite{endo-takemura-list}).  Then we randomly label the vertices to
obtain a decomposable model.  We will discuss random generation of
initial models for $m>8$ in Algorithm 2 below.

In Step 2 we choose the candidate for next decomposable model.
We add or subtract an edge and determine whether the resulting graph
is chordal.  If it is chordal we evaluate its AIC.  For evaluating AIC
we need to obtain the set of cliques and the set of separators.  
Chordality of a graph is determined by obtaining a
perfect elimination scheme and the set of cliques and the 
separators are obtained by ``Maximum cardinality search''
algorithm (\cite{blair-peyton}).  

\bigskip
\indent 
{\bf Algorithm 1} \qquad Model selection of decomposable models.\\
\indent
Input: Microdata $D$, List of non-isomorphic chordal graphs ${\cal L}_m$ with $m$ vertices\\
\indent
Output: Model $M$ with local minimum AIC.\\
\hspace*{1.0cm} $1.$\ Choose a chordal graph from $H\in {\cal L}_m$ at random;\\
\hspace*{1.7cm} \ Label vertices of $H$ at random and obtain a chordal
graph $G_{\rm next}$;\\
\hspace*{1.7cm} \ $A_{\rm next}\leftarrow \text{AIC of } M(G_{\rm next})$;\\
\hspace*{1.0cm} $2.$\ {\bf while} $f=true$ {\bf do}\\
\hspace*{1.7cm} \ $f=false$;\\
\hspace*{1.7cm} \ $G\leftarrow G_{\rm next}$;\\
\hspace*{1.7cm} \ ${\bf for}$\ each\ $e\in E(C_m)$ {\bf do}\\
\hspace*{2.4cm} \ ${\bf if}$ $e\in E_G$ ${\bf then}$\\
\hspace*{3.1cm} \ $G^{\prime }\leftarrow G(V,G(E)\setminus e)$\\
\hspace*{2.4cm} \ ${\bf else}$\\
\hspace*{3.1cm} \ $G^{\prime }\leftarrow G(V,G(E)\cup e)$;\\
\hspace*{2.4cm} \ {\bf if} $G^{\prime }$ is chordal {\bf then}\\
\hspace*{3.1cm} \ $A^{\prime }\leftarrow \text{AIC of } M(G^{\prime })$;\\
\hspace*{3.1cm} \ ${\bf if}$ $A^{\prime }<A_{\rm next}$ ${\bf then}$\\
\hspace*{3.8cm} \ $G_{\rm next}\leftarrow G^{\prime}$;\\
\hspace*{3.8cm} \ $A_{\rm next}\leftarrow A^{\prime}$;\\
\hspace*{3.8cm} \ $f\leftarrow true$;\\
\hspace*{1.0cm} $3.$\ Output $M(G)$;\\

For $m>8$ we can propose the following algorithm to generate an initial
decomposable model to replace Step 1 of Algorithm 1.  Given a chordal
graph $G$ with $m$ vertices, we can obtain a chordal graph $G'$ with
$m+1$ vertices by adding the $m+1$'st vertex and connecting it to a
subset of one clique $C$ of $G$.  Since a chordal graph possesses a
perfect sequence of cliques, the above recursive procedure generates
all chordal graphs. 
The following Algorithm 2 outputs the set of cliques of a random chordal graph.
Note that the probability distribution on random
choices in the algorithm is not specified and the distribution of the output 
is not necessarily the uniform distribution over the set
of chordal graphs with $m$ vertices.

\bigskip
\indent 
{\bf Algorithm 2} \qquad ``Random'' chordal graph with $m$ vertices.\\
\indent
Input: $m$ \\
\indent
Output: Set of cliques of a random chordal graph with $m$ vertices.\\
\hspace*{1.0cm} $1.$\ Initialize $\cal C=\emptyset$;\\
\hspace*{1.0cm} $2.$\ {\bf for } $j \leftarrow 1$ {\bf until } $m$
{\bf do}\\
\hspace*{1.7cm} \ Flip a coin;\\
\hspace*{1.7cm} \ {\bf if} heads {\bf then}\\
\hspace*{2.4cm} \  ${\cal C} \leftarrow {\cal C} \cup \{\{j\}\}$\\
\hspace*{1.7cm} \ {\bf else} choose a member $C\in {\cal C}$
and a subset $C' \subset C$ at random;\\
\hspace*{2.4cm} \ ${\bf if}$ $C=C'$ ${\bf then}$\\
\hspace*{3.1cm} \ $C \leftarrow C\cup \{j\}$\\
\hspace*{2.4cm} \ ${\bf else}$\\
\hspace*{3.1cm} \ ${\cal C} \leftarrow {\cal C} \cup \{C' \cup \{j\} \}$;\\
\hspace*{1.0cm} $3.$\ Output  $\cal C$;\\

\section{Per-record identification risk and estimate of the number of
  population uniques}
\label{sec:risk}

When a good fitting decomposable model is chosen we can estimate the
cell probability of a sample unique cell by MLE (\ref{eq:mle}).
Then a natural estimate of the conditional probability that the sample
unique cell $\bi$  is also a population unique is given as
\begin{equation}
\label{eq:estimate-pop-unique}
(1-\hat p_{{\rm MLE}}(\bi))^{N-n},
\end{equation}
where $N$ is the population size and $n$ is the sample size.  
(\ref{eq:estimate-pop-unique}) is the estimated probability  that none
of the remaining $N-n$ individuals in the population fall into cell $\bi$, 
under the assumption that individuals fall into cells independently
from each other according to the estimated probability distribution.
The number of population uniques in the sample can be estimated as
\[
\sum_{\bi:\text{sample unique}} (1-\hat p_{{\rm MLE}}(\bi))^{N-n} .
\]

In Table \ref{result1} we show two models with smallest values of AIC
by applying Algorithm 1 100 times to our example. Algorithm 1 converged
after a few transitions and it seems to be very practical.  These two
models were also most frequently obtained from Algorithm 1.  In both
models, the separator $\{6\}$ has multiplicity $2$ as indicated by the
repetition in the table.  The estimated numbers of population uniques
(48.867, 40.51) are between those of Ewens model and Pitman model and
seem to be reasonable. The variable 6 (Spouse present/absent) is
contained in many cliques, which can be explained by its high
correlation with other variables and yet small degrees of freedom.  On
the other hand variable 5 (Place of birth) is contained in a single
clique (i.e.\ it is a simplicial vertex), which is also reasonable.

Furthermore the sample uniques with very small estimated cell
probabilities ($p(i) \le 10^{-8}$) are common to these two
models.  We might consider some disclosure control measure for about
20 sample uniques with estimated cell probability less than $10^{-7}$.

\begin{table}[htb]
\caption{Chosen models}
\label{result1}
\begin{center}
\begin{tabular}{|c|c|c|}\hline
 & Model 1& Model 2\\ \hline
Number of times chosen&11 & 7\\ \hline
AIC/2 & 13869.07&13984.97\\ \hline
log likelihood &  $-12141.07$ & $-12013.97$\\ \hline
degrees of freedom & 1728 & 1971\\ \hline
estimated \# of population uniques 
&48.867& 40.515\\ \hline
cliques & \{1,2,6\},\{1,6,7\},\{2,6,8\},&\{1,6,7\},\{3,6,7\},\{1,6,8\},\\
&\{3,6,7\},\{4,6\},\{5,6\} & \{2,8\},\{4,6\},\{5,6\}\\ \hline
separator&\{1,6\},\{2,6\},\{6,7\},\{6\},\{6\}&\{1,6\},\{6,7\},\{6\},\{6\},\{8\}\\ \hline
cell probability estimates& frequencies &  frequencies\\ 
$10^{-2}$ to $10^{-3}$ & 0 &0\\
$10^{-3}$ to $10^{-4}$ & 352 &351\\
$10^{-4}$ to $10^{-5}$ & 1092 &1117\\
$10^{-5}$ to $10^{-6}$ & 599 &600\\
$10^{-6}$ to $10^{-7}$ & 179 &158\\
$10^{-7}$ to $10^{-8}$ & 19 &15\\
$10^{-8}$ to $10^{-9}$ & 2 &2\\
$10^{-9}$ to $10^{-10}$ & 0 &0\\ \hline
\end{tabular}
\end{center}
\end{table}

\section{Swappability of risky records}
\label{sec:swappability}

In Table \ref{result1} two records have the estimated cell probability
of less than $10^{-8}$.  They probably need some disclosure control.
In this paper we propose to swap some observations of these records
with other records of the data set.  Since we have found a
decomposable model with a good fit, it is desirable to swap the
observations such that the marginal frequencies for the cliques of
the chosen model is not disturbed.  In \cite{takemura-hara} we give some
necessary and sufficient conditions for swappability of a particular
sample unique record with some other record without disturbing a given
set of marginals.  

For a decomposable model, a simple method for searching another record
for swapping can be described as follows.  Let $\bi$ be a sample
unique record, such that we want to swap some observations of this
record with another record.  Let $\cal C$ be the set of cliques of a
chosen model and let $\cal S$ denote the set of minimal vertex
separators.  Write each separator $S$ as the intersection of two
cliques $S=C \cap C'$.  We consider all triples $(C,C',S)$ such that 
$S=C \cap C'$.
For example in Model 1 in Table \ref{result1}
all possible  ways of writing separators are as follows.

\allowdisplaybreaks{
\begin{align*}
\{1,6\}&= \{ 1,2,6\} \cap \{ 1,6,7\}, \\
\{2,6\}&=\{ 1,2,6\} \cap \{ 2,6,8\},\\
\{6,7\}&=\{ 1,6,7\} \cap \{ 3,6,7\},\\
\{6\}&=\{ 1,2,6\}\cap \{ 3,6,7\} =  \{ 1,2,6\}\cap \{4,6\} = 
    \{ 1,2,6\}\cap \{5,6\}\\
&=\{ 1,6,7\} \cap \{ 2,6,8\} = \{ 1,6,7\} \cap \{4,6\} = \{ 1,6,7\}
\cap \{5,6\}\\
&= \{ 2,6,8\} \cap \{ 3,6,7\} = \{ 2,6,8\} \cap \{4,6\} = \{ 2,6,8\}
\cap \{5,6\}\\
&= \{ 3,6,7\} \cap \{4,6\} = \{ 3,6,7\}\cap \{5,6\} = \{4,6\}
\cap \{5,6\}.
\end{align*}
}

For a particular sample unique record $\bi$, we search other records
$\bj\neq \bi$ such that for some $(C,C',S)$ we have
\begin{equation}
\label{eq:swappability}
\bi_S=\bj_S, \quad \bi_C \neq \bj_C, \quad \bi_{C'} \neq \bj_{C'}
\end{equation}
If we find some $\bj$ and some $(C,C',S)$ such that
(\ref{eq:swappability}) holds, then we can swap some observations between
$\bi$ and $\bj$.  

We applied this procedure to 50 sample unique records with small
estimated cell probabilities in Table \ref{result1}.  For both models
of Table \ref{result1} this procedure quickly found other records for
swapping for most of 50 records, including the two records with 
the estimated cell probability of less than $10^{-8}$.  Therefore this
procedure seems to work very well in practice.

Note that (\ref{eq:swappability}) is a sufficient condition for swappability
between $\bi$ and $\bj$ for a decomposable model.
For a full statement of necessary and sufficient conditions for
general hierarchical model see Section 3 of \cite{takemura-hara}.

\section{Concluding remarks}
\label{sec:remarks}

In this paper we proposed a systematic strategy for disclosure risk
evaluation and disclosure control of microdata set by fitting
decomposable models.
We have restricted our attention to decomposable models in view of
computational convenience.  Clearly it is desirable to consider other
hierarchical models such as the model containing all two-factor
interaction terms. Simpler hierarchical model might give a better fit
than more complicated decomposable model.  One strategy we can try is
to look for hierarchical models which improves the fit around a locally
best decomposable model.

We have used AIC for evaluating the fit of the model.  Theoretically
AIC is justified for large sample size.  In disclosure control
problems we are dealing with large and sparse tables and from
theoretical viewpoint use of AIC is not justified .  However in
practice it is simple and seems to work reasonably well.  It is of
interest to investigate other methods of model selection for
evaluating the fit of various models.

In microdata sets of official statistics, there are large number of
structural zeros due to various logical relations between key variables.
In principle we should list all the logical relations and specify
structural zeros before fitting a model.  But this is very
cumbersome.  Also the calculation of degrees of freedom of a model
becomes complicated.  It is desirable to develop some practical methods
to deal with structural zeros in some automatic way.

If we want to swap some observations from a sample unique record $\bi$
and if we can find many other records $\bj$ for swapping, it might be
desirable to use $\bj$ which is close to $\bi$ in some sense.  In
\cite{takemura2002} we considered swapping of observations between
close records by introducing an appropriate distance function between
records.

\bigskip
\noindent {\bf Acknowledgment} \quad The approach of this paper was
suggested in a talk by Stephen Fienberg \cite{fienberg-talk} at
University of Tokyo in May 2003 and we are very grateful to his
insights. It took us a long time to implement the whole strategy based on
his suggestions.

\bibliographystyle{plain}
%\bibliographystyle{acmtrans-ims}
% \bibliography{%
%   decomposable%
% %  decomposable0,%
% %  swap,%
% %  groebner%
% }

\end{document}